\documentclass{amsart}
\usepackage[all]{xy}
\usepackage{enumerate}
\usepackage{graphics}
\usepackage[textwidth=11.1cm,textheight=17.1cm]{geometry}
\DeclareMathOperator{\PSL}{PSL}
\DeclareMathOperator{\PGamma}{P\Gamma}
\newcommand\Q{\mathbb{Q}}
\newcommand\N{\mathbb{N}}
\newcommand\R{\mathbb{R}}
\newcommand\C{\mathbb{C}}
\newcommand\Z{\mathbb{Z}}
\newcommand{\sceq}{\mathrel{\mathop:}=}
\DeclareMathOperator{\Ima}{Im}
\DeclareMathOperator{\Rea}{Re}
\newcommand{\textmat}[4]{\left(\begin{smallmatrix} #1&#2 \\ #3&#4
\end{smallmatrix}\right)}
\DeclareMathOperator{\id}{id}
\DeclareMathOperator{\pr}{pr}
\newcommand\wrt{\mbox{w.\,r.\,t.\@ }}
\newcommand\ie{\mbox{i.\,e., }}
\newcommand\resp{\mbox{resp.\@}}
\newcommand{\eps}{\varepsilon}
\DeclareMathOperator{\Cod}{Cod}
\DeclareMathOperator{\NC}{NC}
\DeclareMathOperator{\Ext}{ext}
\DeclareMathOperator{\vb}{vb}
\DeclareMathOperator{\vc}{vc}

\theoremstyle{plain}
\newtheorem{proposition}{Proposition}[section]
\newtheorem{lemma}[proposition]{Lemma}
\newtheorem{theorem}[proposition]{Theorem}

\begin{document}

\title[Symbolic dynamics for the geodesic flow]{Symbolic dynamics for the geodesic flow on locally symmetric orbifolds of rank one}

\author{J. Hilgert and A. D. Pohl}

\address{Institut f\"ur Mathematik, Universit\"at Paderborn,
33095 Paderborn, Germany
}
\email{\{hilgert,pohl\}@math.upb.de}

\begin{abstract}
We present a strategy for a geometric construction of cross sections
for the geodesic flow on locally symmetric orbifolds of rank one. We
work it out in detail for $\Gamma\backslash H$, where $H$ is the
upper half plane and $\Gamma=\PGamma_0(p)$, $p$ prime. Its
associated discrete dynamical system naturally induces a symbolic
dynamics on $\R$. The transfer operator produced from this symbolic
dynamics has a particularly simple structure.
\end{abstract}

\maketitle

\section{Introduction}\label{intro}
We consider the upper half plane $H\sceq \{ z\in \C \mid \Ima z >
0\}$ with the Riemannian metric given by the line element $ds^2 =
y^{-2} (dx^2+dy^2)$ as a model for two-dimensional real hyperbolic
space. The group of orientation-preserving isometries can be
identified with $\PSL(2,\R)$ via the action
\[
g\cdot z = \frac{az+b}{cz+d}
\]
for $g=\textmat{a}{b}{c}{d}\in\Gamma$ and $z\in H$. The geodesic
compactification $\overline H$ of $H$ will be identified with the
one-point compactification of the closure of $H$ in $\C$, hence
$\overline H = \{ z\in\C \mid \Ima z \geq 0\} \cup \{\infty\}$. The
action of $\PSL(2,\R)$ extends continuously to $\partial H = \R\cup
\{\infty\}$.

The geodesics on $H$ are the semicircles centered on the real line
and the vertical lines. All geodesics shall be oriented and
parametrized by arc length. The (unit speed) \textit{geodesic flow}
on $H$ is the dynamical system
\[
\Phi \colon \left\{
\begin{array}{ccc}
\R\times SH & \to & SH
\\
(t,v) & \mapsto & \gamma'_v(t),
\end{array}\right.
\]
where $SH$ denotes the unit tangent bundle (sphere bundle) of $H$
and $\gamma_v\colon\R\to H$ is the unique (unit speed) geodesic that
satisfies $\gamma'_v(0)=v$.

Let $\Gamma$ be a properly discontinuous subgroup of $\PSL(2,\R)$.
The unit tangent bundle $SY$ of the locally symmetric orbifold $Y
\sceq \Gamma\backslash H$ shall be identified with $\Gamma\backslash
SH$, where we consider the induced $\Gamma$-action on $SH$. Let
$\pi\colon H \to Y$ and $\pi\colon SH\to SY$ denote the canonical
projection maps. The geodesics on $Y$ are in bijection with the
$\Gamma$-equivalence classes of geodesics on $H$, and the geodesic
flow on $Y$ is given by
\[
\widehat\Phi \sceq \pi\circ \Phi \circ(\id\times \pi^{-1}) \colon \R\times SY \to SY.
\]
Here, $\pi^{-1}$ shall be read as an arbitrary section of $\pi$. The
definition of $\widehat\Phi$ is independent of this choice. For the
modular group $\PSL(2,\Z)$, Series \cite{Series}
geo\-metrically constructed an amazingly simple cross
section\footnote{The concepts from symbolic dynamics are explained
in Section~\ref{symdyn}.}  for the geodesic flow on the modular
surface $\PSL(2,\Z)\backslash H$. Its associated discrete dynamical
system is naturally related to a symbolic dynamics on $\R$, namely
by continued fraction expansions of the limit points of the
geodesics. The Gauss map is a generating function for the future
part of this symbolic dynamics. In \cite{Mayer_selberg}, Mayer
investigated the transfer operator $\mathcal L_\beta$ with parameter
$\beta$ of the Gauss map. His work and that of Lewis and Zagier
\cite{Lewis_Zagier} have shown that there is an isomorphism
between the space of Maass cusp forms for $\PSL(2,\Z)$ with
eigenvalue $\beta(1-\beta)$ and the space of real-analytic
eigenfunctions of $\mathcal L_\beta$ that have eigenvalue $\pm 1$
and satisfy certain growth conditions.

In this article we present, for the examples $\Gamma=\PGamma_0(p)$,
$p$ prime, a geometric construction of a cross section $\widehat
C_{\mathrm{red}}$ for the geodesic flow on $\Gamma\backslash H$ such
that, as in Series' work, its associated discrete dynamical system
$(\widehat C_{\mathrm{red}}, R)$ naturally gives rise to a symbolic
dynamics on $\R$. More precisely, $(\widehat C_{\mathrm{red}}, R)$
is conjugate to a discrete dynamical system $(\widetilde D,
\widetilde F)$, where $\widetilde D$ is an open subset of
$\R\times\R$. The domain $\widetilde D$ will be seen to naturally
decompose into a finite disjoint union $\bigcup_{\alpha\in A}
\widetilde D_\alpha$ of open subsets $\widetilde D_\alpha$ such that
$\widetilde F\vert_{\widetilde D_\alpha}\colon \widetilde D_\alpha
\to \widetilde F(\widetilde D_\alpha)$ is a diffeomorphism of the
form $(x,y) \mapsto (h_\alpha^{-1}x, h_\alpha^{-1}y)$ for some
$h_\alpha\in\Gamma$. Hence we get a canonical symbolic dynamics for
$\widehat\Phi$ on the alphabet $A$ or $\{h_\alpha \mid \alpha\in
A\}$. If $\pr_x$ denotes the projection onto the first component,
then the sets $D_\alpha \sceq \pr_x(\widetilde D_\alpha)$ are
pairwise disjoint. This means that $\widetilde F$ is effectively
determined by the map $F\colon D\to D$, where $D\sceq
\pr_x(\widetilde D)$ and $F\vert_{D_\alpha} \sceq h_\alpha^{-1}$. In
addition, $F$ is a generating function for the future part of the
symbolic dynamics. Because of the particular structure of $F$, its
associated transfer operator can be described completely by a finite
number of group elements in $\Gamma$.

For the modular surface, this construction gives a two-term transfer
operator from which one can easily read off the three-term
functional equation used in the proof of the Lewis-Zagier
isomorphism mentioned above and whose eigenfunctions with eigenvalue
$1$ that satisfy some growth conditions are in bijection with the
Maass cusp forms.

There are two main methods to construct cross sections and symbolic
dynamics. The \textit{geometric coding} consists in taking a
fundamental domain for $\Gamma$ in $H$ with side pairing and the
sequences of sides cut by a geodesic as coding sequences.  The cross
section is a set of unit tangent vectors based at the boundary of
the fundamental domain. In general, it is very difficult, if not
impossible, to find a conjugate dynamical system on the boundary. In
contrast, the \textit{arithmetic coding} starts with a symbolic
dynamics on (parts of) the boundary which has a generating function
for the future part similar to $F$ from above and asks for a cross
section that reproduces this function. Usually, writing down such a
cross section is a non-trivial task. A good overview of geometric
and arithmetic coding is the survey article
\cite{Katok_Ugarcovici}.

Our method, which we call \textit{cusp expansion}, is sketched in
the following. The description applies to $\Gamma=\PGamma_0(p)$,
$\Gamma=\PSL(2,\Z)$ and some other groups as will be discussed in
Section~\ref{outlook}.  We fix a so-called Ford fundamental domain
$\mathcal F$ for $\Gamma$ in $H$. To each vertex $v$ (that is, an
inner vertex or a representative other than $\infty$ of a cusp) of
$\mathcal F$ we assign a set $A(v)$, which we refer to as a precell.
Two precells overlap at most at boundaries, and the union of all
precells coincides with $\overline{\mathcal F}$. Then we glue
together some $\Gamma$-translates of precells to obtain a cell
$B(v)$ assigned to $v$. The purpose of these cells is to construct a
family of finite-sided $n$-gons with all vertices in $\partial H$
such that each cell has two vertical boundary components (in other
words, $\infty$ is a vertex of each cell) and each non-vertical
boundary component of a cell is a $\Gamma$-translate of a vertical
boundary component of some cell. Further, the family of all cells
shall give a tiling of $H$ in the sense that the $\Gamma$-translates
of all cells cover $H$ and two $\Gamma$-translates of cells are
either disjoint or identical or coincide at exactly one boundary
component. Let $\widetilde C$ denote the set of unit tangent vectors
based at the boundary of some cell but that are not tangent to this
boundary. We will see that $\widehat C \sceq \pi(\widetilde C)$ is a
cross section. In general, $(\widehat C, R)$ is not conjugate to a
discrete dynamical system $(\widetilde D, \widetilde F)$ with
$\widetilde D \subseteq \R\times\R$. But $\widehat C$ contains a
subset $\widehat C_{\mathrm{red}}$ for which $(\widehat
C_{\mathrm{red}}, R)$ is naturally conjugate to a dynamical system
on some subset of $\R\times \R$. The set $\widehat C_{\mathrm{red}}$
is itself a cross section and can be constructed effectively from
$\widehat C$. For the proofs of these properties we extend the
notions of precells and cells to $SH$. Each precell $A(v)$ induces a
precell $\widetilde A(v)$ in $SH$ is a geometric way. From these
precells we construct (in an effective way which, however, involves
choices) a finite family $\{ \widetilde B(v_k)\}_{k\in A}$ of cells
in $SH$ such that each $\widetilde B(v_k)$ projects to a
$\Gamma$-translate of a cell in $H$. The union of all $\widetilde
B(v_k)$ turns out to be a fundamental set for $\Gamma$ in $SH$. This
property and the interplay of the cells in $SH$ and $H$ imply the
properties of $\widehat C$ and $\widehat C_{\mathrm{red}}$. The
decomposition of $\mathcal F$ into precells and the construction of
cells from precells is based on \cite{Vulakh}. Our
constructions differ from Vulakh's in two aspects: one is the way in
which cusps are treated, the other difference is that we extend the
considerations to cells in the unit tangent bundle.

Section~\ref{symdyn} contains the necessary background on symbolic
dynamics. In Section~\ref{cuspexp} we develop the cusp expansion for $\PGamma_0(p)$ in detail and apply it to $\PSL(2,\Z)$. We end, in
Section~\ref{outlook}, with a brief discussion on potential
extensions of this method to other locally symmetric orbifolds of
rank one.

\section{Symbolic dynamics}\label{symdyn}

Let $\widehat C$ be a subset of $SY$. A geodesic $\hat\gamma$ on $Y$
is said to \textit{intersect $\widehat C$ infinitely often in
future} if there is a sequence $(t_n)_{n\in\N}$ with
$\lim_{n\to\infty} t_n = \infty$ and $\hat\gamma'(t_n) \in \widehat
C$ for all $n\in\N$. Analogously, $\hat\gamma$ is said to
\textit{intersect $\widehat C$ infinitely often in past} if we find
a sequence $(t_n)_{n\in\N}$ with $\lim_{n\to\infty} t_n=-\infty$ and
$\hat\gamma'(t_n) \in\widehat C$ for all $n\in\N$. Let $\mu$ be a
measure on the space of geodesics on $Y$. A \textit{cross section}
$\widehat C$ (\wrt $\mu$) for the geodesic flow $\widehat \Phi$ is a
subset of $SY$ such that
\begin{enumerate}[(C1)]
\item \label{C1} $\mu$-almost every geodesic $\hat\gamma$ on $Y$ intersects
$\widehat C$ infinitely often in past and future,
\item \label{C2} each intersection of $\hat\gamma$ and $\widehat C$
is \textit{discrete in time:} if $\hat\gamma'(t) \in\widehat C$,
then there is $\eps > 0$ such that
$\hat\gamma'( (t-\eps, t+\eps) )\cap \widehat C = \{ \hat\gamma'(t) \}$.
\end{enumerate}
Each Borel probability measure on $SY$ which is invariant under $\{
\widehat \Phi_t\}_{t\in\R}$ induces a measure on the space of
geodesics. We shall see that the cross sections constructed in
Section~\ref{cuspexp} are indeed cross sections \wrt any choice of a
measure of this type.

Let $\pr\colon SY\to Y$ denote the canonical projection on base
points. If $\pr(\widehat C)$ is a totally geodesic submanifold of $Y$
and $\widehat C$ does not contain elements tangent to $\pr(\widehat
C)$, then $\widehat C$ automatically satisfies (C\ref{C2}).

Suppose that $\widehat C$ is a cross section for $\widehat \Phi$. If
$\widehat C$ in addition satisfies the property that \textit{each}
geodesic intersecting $\widehat C$ at all intersects it infinitely
often in past and future, then $\widehat C$ will be called a
\textit{strong cross section}, otherwise a \textit{weak cross
section}. Clearly, every weak cross section contains a strong cross
section.

The \textit{first return map} of $\widehat\Phi$ \wrt a strong cross
section $\widehat C$ is the map
\[
R\colon\left\{
\begin{array}{ccc}
\widehat C & \to & \widehat C
\\
\hat v & \mapsto & \widehat{\gamma_v}'(t_0)
\end{array}\right.
\]
where $\pi(v) = \hat v$, $\pi(\gamma_v) = \widehat{\gamma_v}$, the
geodesic $\gamma_v$ is chosen as in the definition of $\Phi$, and
$\widehat{\gamma_v}'(t)\notin \widehat C$ for all $t\in (0,t_0)$. In
this case, $t_0$ is the \textit{first return time} of $\hat v$
resp.\@ of $\widehat{\gamma_v}$. For a weak cross section $\widehat
C$, the first return map can only be defined on a subset of
$\widehat C$. In general, this subset is larger than the maximal
strong cross section contained in $\widehat C$.

Suppose that $\widehat C$ is a strong cross section and let $A$ be
an at most countable set. Decompose $\widehat C$ into a disjoint
union $\bigcup_{\alpha\in A}\widehat C_\alpha$. To each $\hat v\in
\widehat C$ we assign the (two-sided infinite) \textit{coding
sequence} $(a_n)_{n\in\Z} \in A^\Z$ defined by
\[
a_n \sceq \alpha \quad\text{iff $R^n(\hat v) \in \widehat C_\alpha$.}
\]
Note that $R$ is invertible  and let $\Lambda$ be the set of all
sequences that arise in this way. Then $\Lambda$ is invariant under
the left shift $\sigma\colon A^\Z \to A^\Z$, $\left( \sigma(
(a_n)_{n\in\Z} )\right)_k \sceq a_{k+1}$. Suppose that the map
$\widehat C\to \Lambda$  is also injective, which it will be in our
case. Then we have the natural map $\Cod\colon \Lambda \to \widehat
C$ which maps a coding sequence to the element in $\widehat C$ it
was assigned to. Obviously, $R\circ\Cod = \Cod\circ\sigma$. The pair
$(\Lambda,\sigma)$ is called a \textit{symbolic dynamics} for
$\widehat\Phi$. If $\widehat C$ is only a weak cross section and
hence $R$ is only partially defined, then $\Lambda$ also contains
one- or two-sided finite coding sequences.

Let $C'$ be a set of representatives for the cross section $\widehat
C$, that is, $C'\subseteq SH$ and $\pi\vert_{C'}$ is a bijection
$C'\to \widehat C$. For each $v\in C'$ let $\gamma_v$ be the
geodesic determined by $v$. Define $\tau\colon \widehat C\to
\partial H\times
\partial H$ by $\tau(\hat v) \sceq (\gamma_v(\infty),
\gamma_v(-\infty))$ where $v=(\pi\vert_{C'})^{-1}(\hat v)$. For some
cross sections $\widehat C$ it is possible to choose $C'$ in a such
way that $\tau$ is a bijection between $\widehat C$ and some subset
$\widetilde D$ of $\R\times\R$.  In this case the dynamical system
$(\widehat C, R)$ is conjugate to $(\widetilde D, \widetilde F)$,
where $\widetilde F  \sceq \tau\circ R\circ \tau^{-1}$ is the
induced selfmap on $\widetilde D$ (partially defined if $\widehat C$
is only a weak cross section). Moreover, to construct a symbolic
dynamics for $\widehat\Phi$, one can start with a decomposition of
$\widetilde D$ into pairwise disjoint subsets $\widetilde D_\alpha$,
$\alpha\in A$.

Finally, let $(\Lambda, \sigma)$ be a symbolic dynamics with
alphabet $A$. Suppose that we have a map $i\colon\Lambda \to D$ for
some $D\subseteq \R$ such that $i( (a_n)_{n\in\Z})$ depends only on
$(a_n)_{n\in\N_0}$, a (partial) selfmap $F\colon D\to D$, and a
decomposition of $D$ into a disjoint union $\bigcup_{\alpha\in A}
D_\alpha$ such that
\[ F( i( (a_n)_{n\in\Z} )) \in D_\alpha \quad\Leftrightarrow\quad a_1 = \alpha\]
for all $(a_n)_{n\in\Z}\in\Lambda$. Then $F$, more precisely the
triple $(F,i, (D_\alpha)_{\alpha\in A})$, is called a
\textit{generating function for the future part} of
$(\Lambda,\sigma)$. If such a generating function exists, then the
future part of a coding sequence is independent of the past part.

\section{Cusp expansions}\label{cuspexp}

\begin{figure}[h]
\begin{center}
\includegraphics{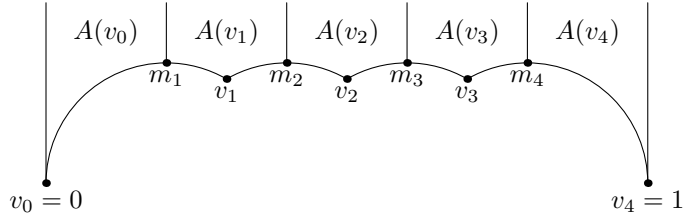}
\end{center}
\caption{Decomposition of $\overline{\mathcal F}$ into precells for $p=5$.}
\label{decomp}
\end{figure}

We consider the groups
\[
\Gamma\sceq \PGamma_0(p) \sceq \left\{ \textmat{a}{b}{c}{d} \in \PSL(2,\Z)
\left\vert\ c\equiv 0\mod p \vphantom{ \textmat{a}{b}{c}{d} }\right.\right\}
\]
for $p$ prime. A subset $\mathcal F$ of $H$ is  a
\textit{fundamental region} for $\Gamma$ in $H$ if $\mathcal F$ is
open, $g\mathcal F\cap \mathcal F = \emptyset$ for all
$g\in\Gamma\setminus\{\id\}$ and
$H=\bigcup_{g\in\Gamma}g\overline{\mathcal F}$. If, in addition,
$\mathcal F$ is connected, it is a \textit{fundamental domain}. A
\textit{fundamental set} for $\Gamma$ in $H$ is a subset of $H$
which contains precisely one representative of each $\Gamma$-orbit.
Clearly, each fundamental region is contained in a fundamental set.
Changing $H$ into $SH$ in these definitions, one gets the notions of
fundamental region, domain and set for $\Gamma$ in $SH$. The
\textit{isometric sphere} of an element
$g=\textmat{a}{b}{c}{d}\in\Gamma$ is the set $I(g) \sceq \{ z\in H
\mid |cz+d|=1\}$, its \textit{exterior} is defined by $\Ext I(g)
\sceq \{ z\in H \mid |cz+d|>1\}$. Note that $I(g)$ is a (complete)
geodesic arc if $g$ does not fix $\infty$. The stabilizer of
$\infty$ in $\Gamma$ is $\Gamma_\infty = \left\{
\textmat{1}{b}{0}{1} \left\vert\ b\in\Z
\vphantom{\textmat{1}{b}{0}{1} }\right.\right\}$. Hence $\mathcal
F_\infty \sceq  \{ z\in H\mid 0 < \Rea z < 1\}$ is a fundamental
domain for $\Gamma_\infty$ in $H$. Following through the proofs of
Thm.~22 in \cite{Ford}, we find that
\[
\mathcal F \sceq \mathcal F_\infty \cap \bigcap_{g\in\Gamma\setminus\Gamma_\infty} \Ext I(g)
\]
is a fundamental domain for $\Gamma$ in $H$. To construct $\mathcal
F$ if suffices to consider the isometric spheres $I_q \sceq \{ z\in
H \mid |pz-q| = 1\}$ for $q=1,\ldots, p-1$. Then $\mathcal F =
\mathcal F_\infty \cap \bigcap_{q=1}^{p-1} \Ext I_q$. From $\mathcal
F$ we can read off that $\Gamma$ has two cusps, represented by
$\infty$ and $v_0\sceq 0$, see Fig.~\ref{decomp}. The points $v_k
\sceq \tfrac{2k+1}{2p} + i\tfrac{\sqrt{3}}{2p}$ for $k=1,\ldots,
p-2$ will be called the \textit{inner vertices} of $\mathcal F$.
Note that $v_k$ is the intersection point of $I_k$ and $I_{k+1}$.
Let $v_{p-1} \sceq 1$ denote the other representative of the cusp
$0$ that is seen by $\mathcal F$. Further let $m_k$ denote the
element of $I_k$ with maximal imaginary part, that is, $m_k\sceq
\tfrac{k}{p} + \tfrac{i}{p}$.

For $a,b\in\overline H$ let $[a,b]$ denote the geodesic arc in $H$
from $a$ to $b$ containing the points $a$ and $b$, let $[a,b)$ be
the geodesic arc from $a$ to $b$ containing $a$ but not $b$, and
define $(a,b]$ and $(a,b)$ analogously. For $a,b\in\R\cup\{\infty\}$
the context will always clarify whether $(a,b)$ refers to a geodesic
arc or a real interval. Then the boundary components of $\mathcal F$
are the geodesic arcs $(v_0,\infty)$, $(v_{p-1},\infty)$ and
$[v_k,v_{k+1}]$ for $k=0,\ldots,p-2$.

For each vertex $v_k$ ($k=0,\ldots, p-1$) the set
\[ A(v_k) \sceq \big\{ z\in \overline{\mathcal F}\ \big\vert\ \tfrac{k}{p} \leq
\Rea z \leq \tfrac{k+1}{p} \big\} \]
is the \textit{precell} attached to $v_k$.
In other words, $A(v_0)$ and $A(v_{p-1})$ are the hyperbolic
triangles with vertices $v_0,m_1,\infty$ resp.\@
$v_{p-1},m_{p-1},\infty$. The precell $A(v_k)$ for an inner vertex
$v_k$ is the hyperbolic quadrangle with vertices
$m_k,v_k,m_{k+1},\infty$. Obviously, the precells overlap only at
boundaries and
\[ \overline{\mathcal F} = \bigcup_{k=0}^{p-1} A(v_k).\]
Let $D$ be a subset of $H$ and $x\in\overline D$. A unit tangent
vector $v$ at $x$ is said to \textit{point into} $D$ if the
geodesic $\gamma$, determined by $\gamma(0)=x$ and $\gamma'(0)=v$
runs into $D$, \ie if there exists $\eps>0$ such that
$\gamma((0,\eps))\subseteq D$. The unit tangent vector $v$ is said
to \textit{point along the boundary} of $D$ if there is $\eps > 0$
such that $\gamma((0,\eps)) \subseteq \partial D$. It is said to
\textit{point out of} $D$ if it points into $H\setminus D$.

For each $A(v_k)$ let $\widetilde A(v_k)$ be the set of unit tangent
vectors that point into its interior $A(v_k)^\circ$. Let $\vb(\widetilde A(v_k))$ be the set of unit tangent vectors based on $\partial A(v_k)$ that point along $\partial A(v_k)$. The set $\vb(\widetilde A(v_k))$ is called the \textit{visual boundary} of $\widetilde A(v_k)$. Further, $\vc(\widetilde A(v_k)) \sceq \widetilde A(v_k) \cup \vb(\widetilde A(v_k))$ is said to be the \textit{visual closure} of $\widetilde A(v_k)$. Then $\vc(\widetilde A(v_k))$ is the set of unit tangent vectors that point into $A(v_k)$.

\begin{lemma}\label{V}
There is a fundamental set $E$ for $\Gamma$ in $SH$ such that
\[
\bigcup_{k=0}^{p-1} \widetilde A(v_k) \subseteq E \subseteq \bigcup_{k=0}^{p-1} \vc(\widetilde A(v_k)).
\]
\end{lemma}

One easily checks that the determining element $g\in\Gamma$ of the
isometric sphere $I(g)$ is unique up to left multiplication by
$\Gamma_\infty$ and that $gI(g)=I(g^{-1})$. This implies that for
each $k\in\{1,\ldots, p-1\}$ there is a unique $l\in\{1,\ldots,
p-1\}$ such that $I_k$ is determined by
\[ g_{kl}\sceq \textmat{l}{ -\frac{1+kl}{p} }{p}{-k} \quad \in \Gamma.\]
E.g., for $k=1$ we have $l=p-1$ and $g_{1,p-1} =
\textmat{p-1}{-1}{p}{-1}$. We now define for each $v_k$,
$k=0,\ldots, p-1$, a set $B(v_k)$, which we call a \textit{cell}, as
follows. For $1\leq k \leq p-2$ we set
\[
 B(v_k)  \sceq \bigcup \{ gA(v_l) \mid l=1,\ldots,p-1,\ g\in\Gamma,\ gv_l=v_k \},
\]
and further $B(v_0)  \sceq  A(v_0) \cup g_{p-1,1}A(v_{p-1})$ resp.\@
$B(v_{p-1})  \sceq A(v_{p-1}) \cup g_{1,p-1}A(v_0)$. The
$\Gamma$-translates of the family $\{ B(v_k)\}_{0\leq k\leq p-1}$
cover $H$, since the $\Gamma$-translates of $\{ A(v_k)\}_{0\leq
k\leq p-1}$ do so. Lemma~\ref{ab} shows that the family of cells
meets the other demands from Sec.~\ref{intro} as well.

\begin{lemma}\label{ab}
For $k=0,\ldots, p-1$, the cell $B(v_k)$ is the hyperbolic triangle
with vertices $\frac{k}{p}$, $\frac{k+1}{p}$ and $\infty$. If
$k\in\{0,p-1\}$, then the non-vertical boundary component of
$B(v_k)$ is
\[
g_{p-1,1}\cdot (v_{p-1},\infty)=\big( v_0, \tfrac{1}{p}\big)
\quad\text{\resp}\quad g_{1,p-1}\cdot (\infty, v_0) = \big(\tfrac{p-1}{p}, v_1\big).
\]
If $1\leq k\leq p-2$, then there are unique $a,b\in\{1,\ldots, p-2\}$ such that
\[
B(v_k) = A(v_k) \cup g_{a,k+1}A(v_a) \cup g_{b+1,k}A(v_b).
\]
The non-vertical boundary component of $B(v_k)$ is
\[
g_{a,k+1}\cdot \big(\tfrac{a+1}{p}, \infty\big) = \big( \tfrac{k}{p},
\tfrac{k+1}{p} \big) = g_{b+1,k}\cdot \big( \infty, \tfrac{b}{p} \big).
\]
Moreover, if $g\partial B(v_k) \cap B(v_l) \not=\emptyset$ for some
$g\in\Gamma$ and $k,l\in\{0,\ldots, p-1\}$, then $g\partial
B(v_k)\cap B(v_l) \subseteq \partial B(v_l)$.
\end{lemma}

For $k=0,\ldots, p-1$ let $C'_k$ denote the set of unit tangent
vectors based on $\tfrac{k}{p}+i\R^+$ that point into
$B(v_k)^\circ$,
\begin{figure}[h]
\begin{center}
\includegraphics{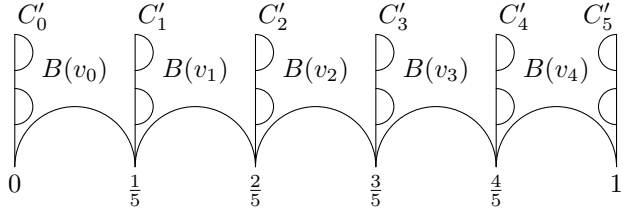}
\end{center}
\caption{The sets $C'_k$ for $p=5$.}
\end{figure}
hence
\[
C'_k = \big\{ X \in SH\ \big\vert\  X = a\tfrac{\partial}{\partial x}
\vert_{\frac{k}{p}+iy} + b\tfrac{\partial}{\partial y}\vert_{\frac{k}{p}+iy},\ a>0,\ b\in\R,\ y>0 \big\}.
\]
Further set
\[
C'_p \sceq \big\{ X \in SH\ \big\vert\  X = a\tfrac{\partial}{\partial x}
\vert_{1+iy} + b\tfrac{\partial}{\partial y}\vert_{1+iy},\ a<0,\ b\in\R,\ y>0 \big\},
\]
which is the set of unit tangent vectors based on $1+i\R^+$ that
point into $B(v_{p-1})^\circ$. Let $\pr\colon SH\to H$ denote the
canonical projection onto base points.

\begin{proposition}
There are pairwise disjoint subsets $\widetilde B(v_k)$ of $SH$, $k=0,\ldots, p$, such that
\begin{enumerate}[(i)]
\item $C'_k\subseteq \widetilde B(v_k)$,
\item the disjoint union $\bigcup_{k=0}^p \widetilde B(v_k)$ is a fundamental set for $\Gamma$ in $SH$,
\item $\overline{\pr(\widetilde B(v_p))} = B(v_{p-1})$ and
$\overline{ \pr(\widetilde B(v_k)) } =B(v_k)$ for $k=0,\ldots, p-1$.
\end{enumerate}
\end{proposition}

\begin{proof}
We pick a fundamental set $E$ for $\Gamma$ in $SH$ of the form as in Lemma~\ref{V}. 
For $z\in\overline{\mathcal F}$ let $E_z$ denote the
set of unit tangent vectors in $E$ that are based at $z$. For each
$k=0,\ldots, p-1$ and each $z\in A(v_k)^\circ$ pick a partition of
$E_z$ into three non-empty disjoint subsets $W^1_{k,z}$,
$W^2_{k,z}$, $W^3_{k,z}$. For $z$ in a non-vertical boundary
component of $A(v_k)$, $k\not= p-1$ and $\Rea(z) \not\in
\{\tfrac{k}{p},\tfrac{k+1}{p}\}$ let $W^1_{k,z} \sceq E_z$ and
$W^2_{k,z} = W^3_{k,z} \sceq \emptyset$. For $z$ contained in a
non-vertical boundary component of $A(v_{p-1})$ and $\Rea(z)
\notin\{ \tfrac{p-1}{p}, 1\}$ divide $E_z$ into two non-empty
disjoint subsets $W^1_{p-1,z}$, $W^3_{p-1,z}$ and let
$W^2_{p-1,z}\sceq \emptyset$. For $z$ with $\Rea(z)=k/p$ for some
$k\in\{0,\ldots, p-1\}$ let $W^1_{k,z}$ be the subset of $E_z$ of
all unit tangent vectors that point into $A(v_k)$ and let
$W^3_{k,z} \sceq E_z\setminus W^1_{k,z}$ and $W^2_{k,z} \sceq
\emptyset$. For $z$ with $\Rea(z)=1$ let $W^1_{p-1,z} \sceq E_z$ and
$W^2_{p-1,z}=W^3_{p-1,z} \sceq \emptyset$. Now we define $\widetilde
B(v_k)$ as follows. For $1\leq k \leq  p-2$ let $a,b$ be as in
Lemma~\ref{ab}. Then
\[
\widetilde B(v_k) \sceq \Big(\bigcup_{z\in A(v_k)} W^1_{k,z}\Big)\ \cup\
\Big(g_{a,k+1}\cdot\bigcup_{z\in A(v_a)} W^2_{a,z}\Big)\ \cup\
\Big(g_{b+1,k}\cdot\bigcup_{z\in A(v_b)} W^3_{b,z}\Big).
\]
For $k\in\{0,p-1,p\}$ we set
\begin{align*}
\widetilde B(v_0) & \sceq \Big(\bigcup_{z\in A(v_0)} W^1_{0,z}\Big)\
\cup\ \Big(g_{p-1,1}\cdot \bigcup_{z\in A(v_{p-1})} W^2_{p-1,z}\Big)
\\
\widetilde B(v_{p-1}) & \sceq \Big(\bigcup_{z\in A(v_{p-1})} W^1_{p-1,z}\Big)\ \cup\
\Big(g_{1,p-1}\cdot \bigcup_{z\in A(v_0)} W^2_{0,z}\Big)
\\
\widetilde B(v_p) & \sceq \Big(\bigcup_{z\in A(v_{p-1})} W^3_{p-1,z}\Big)\
\cup\ \Big(g_{1,p-1}\cdot\bigcup_{z\in A(v_0)}W^3_{0,z}\Big).
\end{align*}
By Lemma~\ref{ab} one sees that the $\widetilde B(v_k)$ satisfy all
the claims of the proposition.
\end{proof}

For reasons that will become clear later (see Prop.~\ref{intcod}) we
shift $\widetilde B(v_p)$ and $C'_p$ one unit to the left, hence let
\[
\widetilde B(v_{-1}) \sceq \textmat{1}{-1}{0}{1}\cdot \widetilde B(v_p)
\quad\text{and}\quad C'_{-1} \sceq \textmat{1}{-1}{0}{1} \cdot C'_p.
\]
Further we set $C' \sceq \bigcup_{k=-1}^{p-1}C'_p$ and
\[
\widehat C \sceq \bigcup_{k=0}^{p} \bigcup_{y>0}
S_{\pi\left(\frac{k}{p}+iy\right)}Y\left\backslash\
\left\{ \pi\left( \pm y^2 \tfrac{\partial}{\partial y}\vert_{\frac{k}{p}+iy}\right)\right\}\right. .
\]
The set $\widehat C$ is the union of all unit tangent vectors in
$SY$ which are based at the boundary of some $\pi(B(v_k))$,
$k=0,\ldots, p-1$, but are not tangent to it. We have the following
lemma which shows that $C'$ is a set of representatives for
$\widehat C$.
\begin{lemma}\label{refset}
$\pi\vert_{C'} \colon C'\to \widehat C$ is a bijection.
\end{lemma}
The next major goal is the following theorem.
\begin{theorem}
$\widehat C$ is a cross section \wrt any measure induced by a Borel
probability measure on $SY$ invariant under
$\{\widehat\Phi_t\}_{t\in\R}$.
\end{theorem}
Let $C \sceq \Gamma\cdot C' = \pi^{-1}(\widehat C)$ and $B\sceq
\Gamma\cdot \bigcup_{k=0}^{p-1} \partial B(v_k) = \pr(C)$. A
geodesic $\hat\gamma$ on $Y$ intersects $\widehat C$ iff some (and
hence any) representative of $\hat\gamma$ intersects $C$. Since $C$
is the set of all unit tangent vectors based in $B$ that are not
tangent to $B$, and since $B$ is totally geodesic, a geodesic
$\gamma$ on $H$ intersects $C$ if and only if it intersects $B$
transversely. Hence $\widehat C$ satisfies (C\ref{C2}).

Suppose that the geodesic $\hat\gamma$ intersects $\widehat C$ in
$\hat\gamma'(t_0)$. By Lemma~\ref{refset} there is a unique geodesic
$\gamma$ on $H$ such that $\gamma'(t_0)\in C'$ and $\pi(
\gamma'(t_0) ) = \hat\gamma'(t_0)$. If
\[
s \sceq \min\{t>t_0 \mid \gamma'(t)\in C \} = \min\{ t>t_0 \mid \gamma(t)\in B\}
\]
exists, then $s$ is the first return time of $\hat\gamma'(t_0)$ and
$R(\hat\gamma'(t_0)) = \pi(\gamma'(s))$. We call $\gamma'(s)$,
resp.\@ $\hat\gamma'(s)$, the \textit{next point of intersection} of
$\gamma$ with $C$, resp.\@ of $\hat\gamma$ with $\widehat C$. It
might happen that $\gamma'(s)\in C'$. Then we say that
$\gamma'(t_0)$ and $\hat\gamma'(t_0)$ are elements with
\textit{interior intersection in future}. The \textit{next point of
exterior intersection} of $\gamma$ and $C$ shall be $\gamma'(r)$,
where
\[
r\sceq \min\{ t>t_0 \mid \exists\, g\in\Gamma\setminus\{\id\}\colon \gamma'(t)\in gC'\}
\]
if $r$ exists. Analogously, we define \textit{previous point of
intersection}, \textit{previous point of exterior intersection} and
elements with \textit{interior intersection in past}.

\begin{figure}[h]
\begin{center}
\includegraphics{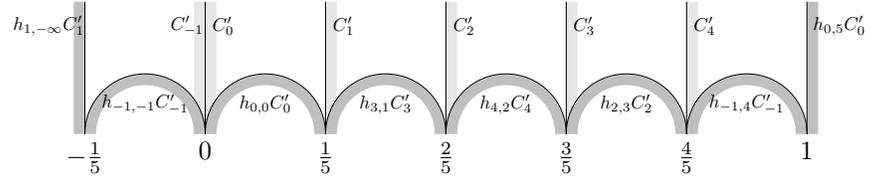}
\end{center}
\caption{The light gray part is the set $C'$, the dark gray parts
are $\Gamma$-translates of the $C'_k$ as indicated. This figure
shows the translates important for determining the next point of
exterior intersection.} \label{nextint}
\end{figure}

\begin{figure}[h]
\begin{center}
\includegraphics{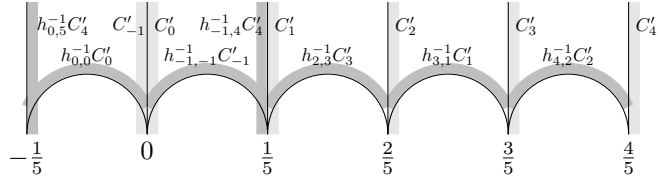}
\end{center}
\caption{Here one can read off the previous point of exterior intersection.}
\label{prevint}
\end{figure}

Prop.~\ref{intcod} below discusses which geodesics intersect
$\widehat C$ at all, which ones intersect it infinitely often, and,
moreover, if there is a next point, resp.\@ was a previous point, of
exterior intersection of the corresponding geodesic on $H$ and on
which copy of $C'$ this point lies. Figs.~\ref{nextint} and
\ref{prevint} show the relevant $\Gamma$-translates of $C'$. To
simplify the statement let $\NC(\Gamma)$ denote the (finite) set of
geodesics on $Y$ which have a representative on $H$ which is tangent
to $B$, and let
\begin{align*}
h_{1,-\infty} & \sceq \textmat{-1}{0}{p}{-1} & h_{c,k} & \sceq g_{c-1,k+1}
\quad\text{for $1\leq k\leq p-2$}
\\
h_{-1,-1} & \sceq \textmat{-1}{0}{p}{-1} & h_{-1,p-1} & \sceq g_{1,p-1}
\\
h_{0,0} & \sceq \textmat{1}{0}{p}{1} & h_{0,p} & \sceq \textmat{1}{1}{0}{1}.
\end{align*}
Finally, let $A\sceq \{-\infty, -1, \ldots, p\}$.

\begin{proposition}\label{intcod} Let $\hat\gamma$ be a geodesic in $Y$.
\begin{enumerate}[{\rm(i)}]
\item \label{propi} $\hat\gamma$ intersects $\widehat C$ iff $\hat\gamma\not\in\NC(\Gamma)$.
\item \label{propii} $\hat\gamma$ intersects $\widehat C$ infinitely often in future iff
$\hat\gamma(\infty)\notin\Gamma\{0,\infty\}$.
\item \label{propiii} $\hat\gamma$ intersects $\widehat C$ infinitely often in past iff
$\hat\gamma(-\infty)\notin\Gamma\{0,\infty\}$.
\item \label{propiv} Suppose that $\hat\gamma$ intersects $\widehat C$ in $\hat\gamma'(t_0)$.
Let $\gamma$ be the unique lifted geodesic on $H$ such that $\gamma'(t_0)\in C'$ and
$\pi(\gamma'(t_0)) = \hat\gamma'(t_0)$.
\begin{enumerate}[{\rm(a)}]
\item \label{propiva} There is a next point of exterior intersection of $\gamma$ and
$C$ iff $\gamma(\infty)\notin\big\{ \frac{k}{p}\ \big\vert\ k=-1,\ldots,p \big\}$.
\item \label{propivb} If $\frac{k}{p} < \gamma(\infty) < h_{c,k}\infty$,
$k\in A\setminus\{-1, p-1\}$, resp.\@ $h_{-1,-1}\infty < \gamma(\infty) < 0$ or
$h_{-1,p-1}\infty <\gamma(\infty)<1$, then the next point of exterior intersection is on
$h_{c,k}\cdot C'_c$.
\item \label{propivc} There was a previous point of exterior intersection of $\gamma$ and $C$
iff $\gamma(\infty) < 0$ and $\gamma(-\infty) \not=\tfrac{1}{p}$, or $\gamma(\infty) > 0$ and
$\gamma(-\infty) \notin \big\{ \frac{k}{p}\ \big\vert\ k=-1,\ldots, p-2\big\}$.
\item \label{propivd} We have the following cases:\\
\begin{tabular}{|c|c|l|}
\hline
$\gamma(-\infty) \in $ & $\gamma(\infty) \in $ & prev pt of ext intsec on
\\\hline
$\big(-\infty, -\tfrac{1}{p}\big)\vphantom{\left(\tfrac1p\right)}$ & $(0,\infty)$ &
\hspace{0.5cm} $h_{0,p}^{-1}\cdot C'_{p-1}$
\\
$\big(-\tfrac{1}{p}, 0\big)$ & $(0,\infty)$ & \hspace{0.5cm} $h_{0,0}^{-1}\cdot C'_0$
\\
$\big(0,\tfrac{1}{p}\big)$ & $\big(\tfrac{1}{p},\infty\big)$ & \hspace{0.5cm} $h_{-1,-1}^{-1}\cdot C'_{-1}$
\\
$\big(0,\tfrac{1}{p}\big)$ & $(-\infty,0)$ & $\hspace{0.5cm} h^{-1}_{-1,-1}\cdot C'_{-1}$
\\
$\big(\tfrac{1}{p},\infty\big)$ & $(-\infty, 0)$ & \hspace{0.5cm} $h_{-1,p-1}^{-1}\cdot C'_{p-1}$
\\
$\big(\tfrac{k}{p}, \tfrac{k+1}{p}\big)$, $1\leq k \leq p-2$ & $\big(\tfrac{k+1}{p},\infty\big)$ &
\hspace{0.5cm} $h_{k+1,b}^{-1}\cdot C'_b$
\\ \hline
\end{tabular}
\end{enumerate}
\end{enumerate}
\end{proposition}

Prop.~\ref{intcod} shows that exactly the geodesics on $Y$ with one
or two limit points in the cusps do not intersect $\widehat C$
infinitely often. Let $E$ denote the set of unit tangent vectors to
these geodesics. Adapting the proof in
\cite{Katok_Ugarcovici_arith}, p.~59, to our situation, we get
that $E$ is a null set \wrt any Borel probability measure on $SY$
which is invariant under $\{\widehat \Phi_t\}_{t\in\R}$. Therefore
$\widehat C$ also satisfies (C\ref{C1}) and hence is a cross
section. The maximal strong cross section contained in $\widehat C$
is $\widehat C_{\mathrm{s}} \sceq \widehat C\setminus E$.

Recall the map $\tau$ from Section~\ref{symdyn}. If $\hat v\in
\widehat C$ is an element with interior intersection in future, then
$\tau(\hat v) = \tau( R(\hat v))$ since the geodesics on $H$
corresponding to $\hat v$ resp.\@ to $R(\hat v)$ only differ by
their parametrization. Likewise, if $\hat v \in\widehat C$ is an
element with interior intersection in past, then $\tau(\hat v) =
\tau(R^{-1}(\hat v))$. Figs.~\ref{nextint} and \ref{prevint} clearly
show that $\widehat C$ contains both, elements with interior
intersection in future and in past. Hence it is impossible to find a
dynamical system on $\partial H$ that is conjugate to $(\widehat C,
R)$. To overcome this problem we reduce $\widehat C$ to another
cross section. Let $P$ denote the set of elements in $\widehat C$
that have interior intersection in future or past and set
\[
\widehat C_{\mathrm{red}} \sceq \widehat C\setminus P\quad\text{resp.}\quad
\widehat C_{\mathrm{red,s}}\sceq \widehat C_{\mathrm{s}}\setminus P.
\]
Prop.~\ref{intcod} implies that each geodesic on $Y$ that intersects
$\widehat C$ infinitely often, also intersects $\widehat
C_{\mathrm{red,s}}$ infinitely often. Therefore both $\widehat
C_{\mathrm{red}}$ and $\widehat C_{\mathrm{red,s}}$ are cross
sections and $\widehat C_{\mathrm{red,s}}$ is a strong cross
section. For $0\leq k \leq p-1$ let
\[
\widetilde D_k \sceq \big(\tfrac{k}{p}, \tfrac{k+1}{p}\big) \times \big(-\infty, \tfrac{k}{p}\big).
\]
Further set
\begin{align*}
\widetilde D_{-\infty} &\sceq \big(-\infty, -\tfrac1p\big) \times (0,\infty),
&\widetilde D_{-1} &\sceq \big(-\tfrac1p,0\big) \times (0,\infty),
\\
\widetilde D_p &\sceq (1,\infty)\times (-\infty,1), & \widetilde D_{k,s}
&\sceq \widetilde D_k \setminus \big(\Gamma\{0,\infty\}\times\Gamma\{0,\infty\} \big)
\end{align*}
and $h_k\sceq h_{c,k}$ for $k\in A$. Define the selfmap $\widetilde
F$ on  $\widetilde D \sceq \bigcup_{k\in A} \widetilde D_k$ resp.\@
$\widetilde F_s$ on $\widetilde D_s \sceq\bigcup_{k\in A} \widetilde
D_{k,s}$ by
\[
\widetilde F\vert_{D_k} \sceq h_k^{-1} \quad\text{resp.}\quad
\widetilde F_s\vert_{D_{k,s}} \sceq h_k^{-1}.
\]
Then Prop.~\ref{intcod} implies that the two diagrams
\[
\xymatrix{ \widehat C_{\mathrm{red}} \ar[r]^R  \ar[d]_{\tau} &
\widehat C_{\mathrm{red}} \ar[d]^{\tau} &&& \widehat
C_{\mathrm{red,s}} \ar[r]^R  \ar[d]_{\tau} & \widehat
C_{\mathrm{red,s}} \ar[d]^{\tau}
\\
\widetilde D \ar[r]^{\widetilde F} & \widetilde D &&& \widetilde D_s\ar[r]^{\widetilde F_s} &
\widetilde D_s
}
\]
commute. The sets $\widetilde D_k$ resp.\@ $\widetilde D_{k,s}$
satisfy all the properties announced in Sec.~\ref{intro}. Therefore
we naturally get a symbolic dynamics with alphabet $A$ and a natural
generating function for its future part.

\subsection{Associated transfer operators}

We restrict ourselves to the strong cross section $\widehat
C_{\mathrm{red, s}}$. Let
\[D  \sceq \pr_x(\widetilde D_s) = \R\setminus\Gamma\{0,\infty\}\] and
$D_k \sceq \pr_x(\widetilde D_{k,s})$ for $k\in A$. The discrete
dynamical system $(\widetilde D_s, \widetilde F_s)$ induces a
symbolic dynamics with generating function $F\colon D\to D$ given by
$F\vert_{D_k} \sceq h_k^{-1}$ for $k\in A$. Its local inverses  are
\[
F_k \sceq \left( F\vert_{D_k}\right)^{-1} \colon F(D_k) \to D_k,\  x\mapsto h_kx.
\]
Then the transfer operator with parameter $\beta$
\[
\mathcal L_\beta \colon \{\text{densities on $D$}\} \to \{\text{densities on $D$}\}
\]
associated to $F$ is given by
\[
(\mathcal L_\beta\varphi)(x) = \sum_{k\in A} F'_k(x)^\beta \varphi( F_k(x) )
\chi_{F(D_k)}(x),
\]
where $\chi_{F(D_k)}$ is the characteristic function of $F(D_k)$,
and the maps $F_k$ and $F'_k$ are extended arbitrarily on
$D\setminus F(D_k)$.

\subsection{The case of the modular surface}

We apply the cusp expansion to the modular group
$\Gamma=\PSL(2,\Z)$. All proofs are omitted as they are analogous to
those for $\PGamma_0(p)$. The set
\[
\mathcal F \sceq \{ z\in H \mid 0<\Rea z < 1,\ |z|>1,\ |z-1|>1\}
\]
is a fundamental domain for $\Gamma$ in $H$.
Note that
\[
\mathcal F = \mathcal F_\infty \cap \bigcap_{g\in\Gamma\setminus\Gamma_\infty} \Ext I(g)
\]
with $\mathcal F_\infty \sceq \{ z\in H\mid 0 < \Rea z < 1\}$. The
point $v\sceq \tfrac12(1+ i\sqrt{3})$ is the only (inner) vertex of
$\mathcal F$. Its associated precell $A(v)$ coincides with
$\overline{\mathcal F}$. The cell $B(v)$ is the hyperbolic triangle
with vertices $0$, $1$, $\infty$. As before, set $Y\sceq
\Gamma\backslash H$ and let $\pi\colon SH\to SY$ denote the
canonical projection. Further, let $\widehat C$ be the set of unit
tangent vectors in $SY$ which are based on $\pi(B(v))$, but that are
not tangent to it. Then $\widehat C$ is a weak cross section with
set of representatives
\[
C' \sceq \big\{ X\in SH\ \big\vert\ X=
a\tfrac{\partial}{\partial x}\vert_{iy} + b\tfrac{\partial}{\partial y}\vert_{iy},\ y>0,\
a>0,\ b\in\R \big\}.
\]
The elements of $C'$ are the unit tangent vectors based on $i\R^+$
that point into $B(v)^\circ$. One easily sees that $\widehat C$
does not contain elements with interior intersection. Let $\widehat
C_s$ denote the maximal strong cross section contained in $\widehat
C$. Then $(\widehat C_s, R)$ is conjugate to $(\widetilde D_s,
\widetilde F_s)$, where $\widetilde D_s \sceq \widetilde
D_{0,s}\cup\widetilde D_{1,s}$ and
\begin{align*}
\widetilde D_{0,s} & \sceq \big((0,1)\setminus\Q\big) \times \big((-\infty,0)\setminus\Q\big), &
\widetilde F_s\vert_{\widetilde D_{0,s}} & \sceq \textmat{1}{0}{-1}{1},
\\
\widetilde D_{1,s} & \sceq \big((1,\infty)\setminus\Q\big) \times \big((-\infty,0)\setminus\Q\big), &
\widetilde F_s\vert_{\widetilde D_{1,s}} & \sceq \textmat{1}{-1}{0}{1}.
\end{align*}
Note that $\Q=\Gamma\infty$. Set $D\sceq \R^+\setminus\Q$, $D_0\sceq
(0,1)\setminus\Q$ and $D_1\sceq (1,\infty)\setminus\Q$. Then the
function $F\colon D\to D$, given by $F\vert_{D_0} \sceq
\textmat{1}{0}{-1}{1}$, $F\vert_{D_1}\sceq \textmat{1}{-1}{0}{1}$,
is the first component of $\widetilde F_s$ restricted to the first
variable. Its associated transfer operator with parameter $\beta$ is
\[
\left(\mathcal L_\beta \varphi\right)(x) = \varphi(x+1) +
(x+1)^{-2\beta}\varphi\big(\tfrac{x}{x+1}\big),\quad x\in D.
\]
If we extend $\mathcal L_\beta$ to act on densities on $\R^+$ by the
obvious formula, then the eigenfunctions of $\mathcal L_\beta$ with
eigenvalue $1$ are precisely the solutions of
\[
 \varphi(x) = \varphi(x+1) + (x+1)^{-2\beta}\varphi\big(\tfrac{x}{x+1}\big),\quad x\in \R^+.
\]
This is exactly the functional equation used in
\cite{Lewis_Zagier}. We do not think that it is a coincidence
that eigenfunctions of Mayer's transfer operator are related to
those of $\mathcal L_\beta$, since $F$ is closely related to the
Farey process or slow continued fractions, whereas the Gauss map is
the generating function for (ordinary) continued fractions.

\section{Outlook}\label{outlook}

The construction of a symbolic dynamics for the geodesic flow on
$\Gamma\backslash H$ via the method of cusp expansion seems to work
for a wide class of subgroups $\Gamma$ of $\PSL(2,\R)$. Some
necessary properties of such a group $\Gamma$ are that $\infty$ is a
cusp and that the boundary of
\[
\bigcap_{g\in \Gamma\setminus\Gamma_\infty}\Ext I(g)
\]
contains the maxima of the relevant isometric spheres. For example,
cusp expansion is applicable to Hecke triangle groups, thus also to
certain non-arithmetic groups. The precells in $H$ can be defined
independently of the choice of a fundamental domain for $\Gamma$ in
$H$. But if there is no close relation between the union of precells
and a fundamental domain for $\Gamma$, we do not expect to be able
to construct symbolic dynamics with this method. On the positive
side, we do expect it to work for certain groups acting on
higher-dimensional real hyperbolic spaces and, more generally, on
rank one symmetric spaces of non-compact type.

\bibliographystyle{amsplain}
\bibliography{hp}

\end{document}